\documentclass{article}
\usepackage[a4paper, total={6in, 8in}]{geometry}
\usepackage{amssymb} % For \mathbb{N} and \mathbb{R}
\usepackage{graphicx} % Required for inserting images
\usepackage{amsmath}
\usepackage{amssymb}
\usepackage{mathrsfs}
\usepackage{cite}
\usepackage{verbatim}
\usepackage{float}
\usepackage{graphicx}
\usepackage{amsthm}
\usepackage{textcomp}
\usepackage{subfig}
\usepackage{esint}
\usepackage{enumerate}
\usepackage{hyperref}
\usepackage{xcolor}

\newcommand{\eps}{\varepsilon} 
\newcommand{\Om}{\Omega} 
\newcommand{\R}{\mathbb{R}} 
\newcommand{\N}{\mathbb{N}}
\newcommand{\mres}{\mathbin{\vrule height 1.6ex depth 0pt width 0.13ex\vrule height 0.13ex depth 0pt width 1.3ex}}
\newtheorem{theorem}{Theorem}
\newtheorem{lemma}{Lemma}
\newtheorem{proposition}{Proposition}
\newtheorem{corollary}{Corollary}
\newtheorem{remark}{Remark}
\newtheorem{definition}{Definition}

\title{On the first eigenvalue of the infinity Fractional Bilaplacian}
	\author{S. Carano \thanks{Dipartimento di Matematica,
Via Sommarive, 5 - 38123 Povo, Trento, Italy.\\
e-mail:  simone.carano@unitn.it} \qquad A. Carbotti \thanks{Dipartimento di Matematica
	e Fisica ``E. De Giorgi'', Universit\`a del Salento,
	Via Per Arnesano, 73100 Lecce, Italy.}}
\date{\today}

\begin{document}

\maketitle

\begin{abstract}
    We study the first eigenvalue of the infinity fractional bilaplacian under hinged and clamped exterior conditions. On bounded Lipschitz domains satisfying a uniform exterior ball condition, we establish existence of minimizers for both problems. In the hinged case, we prove uniqueness up to multiplication by a nonzero constant and identify the minimizers with the fractional torsion function. In the clamped case, an $L^p$-approximation yields a minimizer satisfying a limiting fractional PDE-system with measure data. We also prove a sharp Faber–Krahn inequality for the hinged eigenvalue, and use its explicit expression to obtain double sided estimates for the clamped eigenvalue on a ball.
\end{abstract}


\section{Introduction}
Variational problems in $L^\infty$ have received considerable attention in recent years, due to their connections with nonlinear partial differential equations, optimal Lipschitz extensions, and  game theory. Among these, $L^\infty$-eigenvalue problems represent a particularly interesting class, whose study was initiated by Juutinen, Lindqvist, and Manfredi in their seminal paper \cite{JLM}, where they introduced the $\infty$-eigenvalue problem as the limiting regime of 
\[
\lambda_p(\Omega) = \inf_{u \in W_0^{1,p}(\Omega) \setminus \{0\}} \frac{\|\nabla u\|_{L^p(\Omega)}}{\|u\|_{L^p(\Omega)}}
\]
as $p \to \infty$. Since then, the theory has been developed in several directions, including higher order operators, such as the $\infty$-bilaplacian studied by Katzourakis and Parini in \cite{KP}, as well as extensions to the nonlinear and fractional setting (see, e.g., \cite{CD,BPS,BCF,LL}).

The aim of the present paper is to contribute to this direction by introducing and studying fractional counterparts of the hinged and clamped $L^\infty$-eigenvalue problems in \cite{KP}.
More precisely, for $s\in(0,1)$, $n\in\mathbb{N}$, and an open bounded set $\Omega\subset\mathbb{R}^{n}$, we consider two natural variational formulations associated with the fractional Laplacian of order $2s$:
\begin{equation}
\label{eq:lambda}
\Lambda_{s,\infty}^{H}(\Omega) := \inf_{u \in \mathscr{L}_{0}^{2s,\infty}(\Omega)\setminus\{0\}} \frac{\|(-\Delta)^{s}u\|_{L^{\infty}(\Omega)}}{\|u\|_{L^{\infty}(\mathbb{R}^{n})}}
\end{equation}
and
\begin{equation}
\label{eq:lambda'}
\Lambda_{s,\infty}^{C}(\Omega) := \inf_{u \in \mathscr{W}_{0}^{2s,\infty}(\Omega)\setminus\{0\}} \frac{\|(-\Delta)^{s}u\|_{L^{\infty}(\mathbb{R}^{n})}}{\|u\|_{L^{\infty}(\mathbb{R}^{n})}},
\end{equation}
where
\begin{align*}
&\mathscr L^{2s,\infty}_0(\Om):=\{u\in L^\infty(\R^n): u=0\,\mbox{ in }\,\R^n\setminus\Om, (-\Delta)^s u\in L^\infty(\Om)\},\\
&\mathscr W^{2s,\infty}_0(\Om):=\{u\in L^\infty(\R^n): u=0\,\mbox{ in }\,\R^n\setminus\Om, (-\Delta)^s u\in L^\infty(\R^n)\},
\end{align*}
and $(-\Delta)^s u$ denotes the distributional $s$-fractional Laplacian of $u$. We postpone to Section \ref{sec:not} for all the missing definitions.
The fractional eigenvalues $\Lambda_{s,\infty}^{H}(\Omega)$ and $\Lambda_{s,\infty}^{C}(\Omega)$ serve as the nonlocal generalizations of $\Lambda_{\infty}^{H}(\Omega)$ and $\Lambda_{\infty}^{C}(\Omega)$ introduced in \cite{KP}, respectively the hinged and clamped first eigenvalue of the $\infty$-Bilaplacian. 
We point out that, differently from \cite{KP}, we do not pursue a formal derivation of the fractional infinity Bilaplacian operator. Such a derivation is already highly nontrivial at a formal level and, even if achievable, would not necessarily provide a meaningful characterization of solutions to the associated eigenvalue problem. Indeed, as also occurs in the classical setting in \cite{KP}, the lack of sufficient regularity of minimizers of $L^\infty$-variational problems poses a serious obstacle to interpreting the resulting operator and its associated eigenvalue equation.

Before presenting our main results, we briefly outline the mathematical context of our problem. We consider an $L^\infty$-variational problem involving fractional differential operators of order $2s$. While the theory for integer order operators is relatively well developed, several challenges remain. In the first order case, most of the main difficulties have been successfully addressed (see \cite{A} for Aronsson's seminal work and \cite{Kbook} for a survey). By contrast, the study of higher order problems, initiated more recently in \cite{KM,KP2}, is still at an early stage, as many first order techniques do not extend to this setting. For further references on higher order problems, we refer to \cite{DK,KM2,KM3,CKM}.

The fractional order setting, on the other hand, remains largely unexplored. In \cite{CLM} (see also \cite{BCF}), the authors introduce a notion of fractional infinity Laplacian, extending the Aronsson equation to the fractional framework. This equation arises as the limit, as $p\to\infty$, of the Euler–Lagrange equations associated with the minimization of the $W^{s,p}$-norm. More generally, the asymptotic behavior of functionals on $W^{s,p}$ as $p\to\infty$ has been extensively investigated (see, e.g., \cite{BPS,AMRT,LL}). However, the PDEs obtained in these studies do not appear to be directly connected to fractional order $L^\infty$-variational problems of Dirichlet or eigenvalue type. Motivated by the approach in \cite{CM}, where an $L^\infty$-Dirichlet problem for the fractional Laplacian is studied, we consider the corresponding eigenvalue problems \eqref{eq:lambda} and \eqref{eq:lambda'}, with the aim of deriving a fractional PDE that provides a necessary, and possibly sufficient, condition for minimality. Furthermore, some interesting work in the field of fractional eigenvalue problems, relevant to the content of this paper, appears in \cite{BBB, BLP, BP, BrPaSq}.

Our first main result is on the existence, uniqueness up to a multiplicative constant, and characterization of minimizers for $\Lambda^H_{s,\infty}(\Om)$.
\begin{theorem}\label{thm1}
Let $\Om\subset\R^n$ be a bounded open set with Lipschitz boundary, satisfying a uniform exterior ball condition. Then the eigenvalue problem for  $\Lambda^H_{s,\infty}(\Om)$, defined in \eqref{eq:lambda}, is non-trivial, i.e. $$\Lambda^H_{s,\infty}(\Om)>0$$ and it admits a solution $u_\infty^H\in \mathscr L^{2s,\infty}_0(\Om)$, which is uniquely determined up to a non-zero multiplicative constant. Moreover, $u^H_\infty$ is a multiple of the solution to the Dirichlet problem
\begin{equation}\label{talenti}
\left\{
\begin{aligned}
(-\Delta)^s v&=1\qquad&&\mbox{in }\Om,\\
v&=0&&\mbox{in }\R^n\setminus\Om.
\end{aligned}
\right.
\end{equation}
\end{theorem}

As a consequence of the results in \cite{RS}, the optimal regularity of $u_\infty^H$ is $C^s(\R^n)$, which coincides with the regularity of the competitors in $\mathscr L^{2s,\infty}_0(\Om)$ (as established in Section \ref{reg comp}). The absence of any regularity gain for the minimizer is a characteristic feature of $L^\infty$ variational problems, even in the local setting (compare e.g. \cite{CKM,KM,KP,KP2}).
On the other hand, as shown in Section \ref{exist sec}, the class of competitors for the clamped problem inherits Sobolev regularity in the space $W^{2s,p}(\R^n)$ for all $p\in[2,\infty)$.
This leads to a strict inequality, in general, in the natural comparison
$$
\Lambda_{s,\infty}^{H}(\Omega) \leq \Lambda_{s,\infty}^{C}(\Omega).
$$ 
In fact, when $s>\frac12$, the Morrey-Sobolev embedding yields $$W^{2s,p}(\R^n)\subset C^{1,\alpha}(\R^n),$$ with $\alpha=2s-1-\frac np>0$ (for $p$ large enough). Consequently, every competitor $u\in\mathscr W^{2s,\infty}_0(\Om)$ is globally of class $C^{1,\alpha}$.
In particular, its gradient is well defined on $\partial\Om$, so that the classical clamped condition  $$\partial_\nu u=0\qquad\mbox{ on }\partial\Om$$ is recovered.
As explicit examples show, e.g. for $\Om=B_r$, the solution to \eqref{talenti} does not possess a well defined normal derivative on $\partial\Om$. Moreover, its $s$-Laplacian does not remain bounded outside $\Om$. Thus, none of the competitors for the clamped problem in $B_r$ can be a minimizer for the hinged one.

In our second main result, we state the existence of a minimizer for $\Lambda^C_{s,\infty}(\Om)$, which in turn must satisfy a (fractional) PDE system involving measure data.

\begin{theorem}\label{thm2}
Let $\Om\subset\R^n$ be a bounded open set with Lipschitz boundary, satisfying a uniform exterior ball condition. Then the eigenvalue problem for  $\Lambda^C_{s,\infty}(\Om)$, defined in \eqref{eq:lambda'}, admits a solution $u^C_\infty\in \mathscr W^{2s,\infty}_0(\Om)$. 
Moreover, there exist two measures $\mu,\nu\in \mathcal{M}(\R^n)\setminus\{0\}$, with $\mu\mres \Om=f_\infty \mathscr{L}^n\mres\Om$, for a function $f_\infty\in L^1(\Om)\setminus\{0\}$, such that
\begin{equation}\label{PDE clamped}
\left\{
\begin{aligned}
(-\Delta)^s u^C_\infty&=\Lambda^C_{s,\infty}(\Om)\mbox{sgn}f_\infty\quad&&\mbox{ a.e. in }\Om\setminus\{f_\infty=0\},\\
(-\Delta)^s\mu&=\nu&&\mbox{in }\Om.
\end{aligned}
\right.
\end{equation}
The last equation is understood in the sense of distributions.\\
More generally, we have
\begin{align}\label{EL out Om}
    (-\Delta)^su_\infty^C=\Lambda^C_{s,\infty}(\Om)\frac{d\mu}{d|\mu|}\qquad\mbox{in }\mathrm{supp(\mu)}\setminus\partial\Om.
\end{align}
Finally, the measures $\mu$ and $\nu$ have compact support and finite total variation. Precisely, we have
\begin{align*}
\mbox{supp}(\mu)\Subset\R^n,\qquad |\mu|(\R^n)\leq\frac{1}{\Lambda^C_{s,\infty}(\Om)},
\end{align*}
and
\begin{align*}
&\mbox{supp}(\nu)\subseteq\left\{x\in\overline\Om: |u^C_\infty(x)|=\|u_\infty^C\|_{L^\infty(\Om)}\right\}\Subset\Om,\\
&|\nu|(\R^n)=|\nu|(\Om)\leq 1.
\end{align*}

\end{theorem}
We note that, unlike in the corresponding Dirichlet problem studied in \cite{CM}, we cannot rule out the possibility that the zero set of $f_\infty$ has positive measure. Nevertheless, we expect (and conjecture) that for domains $\Om$ with suitable symmetry properties (e.g. $\Om=B_r$), the zero nodal set of $f_\infty$ is negligible. This would yield the eikonal property
$|(-\Delta)^s u_\infty^C|=\Lambda^C_{s,\infty}(\Om)$ a.e. in $\Om$.\\
 The existence of a minimizer $u_\infty^C$ and the derivation of the fractional PDE system are proved via $L^p$-approximation (i.e. Gamma Convergence). The main idea is to consider and solve the corresponding $p$-version of problem \eqref{eq:lambda'}, and retain the information in the limit as $p\to\infty$. The proof is essentially technical and relies on a combination of delicate regularity estimates for fractional partial differential equations. In particular, it makes extensive use of the global Calderón–Zygmund estimates in \cite{AFLY}, together with regularity results for fractional PDEs in the distributional setting which are developed in Section \ref{prelim}. 
A further delicate point is to ensure that the measures $\mu$ and $\nu$, as well as the function $f_\infty$, are nontrivial. This requires exploiting specific properties of $s$-harmonic functions. More precisely, we rely on a suitable version of the well known result by Dipierro, Savin, and Valdinoci \cite{DSV}, according to which ``all functions are locally $s$-harmonic up to a small error". A recent variant of this principle, established in \cite{CM}, plays a crucial role in our argument.

Finally, we establish a Faber-Krahn inequality for hinged external conditions.

\begin{theorem}\label{FK hinged}
    Let $\Omega$ be a bounded Lipschitz domain satisfying a uniform exterior ball condition. Let $R>0$ be such that $|\Omega|=|B_R|$. Then, one has
    \begin{equation}
    \label{eq:FKH}
    \Lambda^H_{s,\infty}(\Omega)\ge \Lambda^H_{s,\infty}(B_R),
    \end{equation}
    and the equality holds in \eqref{eq:FKH} if, and only if, $\Omega$ is a ball.
\end{theorem}

This result relies heavily on the symmetrization techniques developed in \cite{FV} and can be regarded as a fractional counterpart of Proposition 3.1 in \cite{KP}. In particular, it extends the corresponding Faber–Krahn inequality to the fractional setting, though the lack of the corresponding Talenti inequality does not allow to accomplish to the rigidity of the equality case.
To overcome this difficulty, we recall in Lemma \ref{lem:heat_symmetrization} that $-(-\Delta)^s$ endowed with Dirichlet external conditions is the infinitesimal generator of a rotationally symmetric $2s$-stable Lévy process, and we use the representation of the torsion functions in terms of the fractional Dirichlet heat flow in Lemma \ref{lem:smalltimeasymp}, reducing this way to a comparison of the initial rates of the heat loss in the flows through a set and its symmetral.

It is worth pointing out that the situation is considerably more delicate for fractional clamped boundary conditions. Indeed, unlike in the local setting, we are currently unable to establish the corresponding Faber–Krahn inequality in the fractional framework. Nevertheless, the explicit expression of the hinged eigenvalue of the ball allows to prove the following estimate.
\begin{proposition}
\label{prop:estimatesforclamp}
    Let $\Lambda_{s,\infty}^C(B_R)$ be the clamped eigenvalue for the ball $B_R$. Then the following estimates hold true
    \begin{equation}
    \label{eq:estimatesclampedeigenvalue}
\Lambda_{s,\infty}^H(B_R)\le\Lambda_{s,\infty}^C(B_R)\le(1+s)\max\left\{1,\frac{2s}{n}\right\}\Lambda_{s,\infty}^H(B_R).
    \end{equation}
\end{proposition}

The rest of the paper is organized as follows. In Section 2, we fix notation and function spaces. Section 3 contains preliminary results on the regularity of distributional solutions to the fractional Dirichlet problem. Sections 4 and 5 are devoted to the proofs of Theorem 1 and Theorem 2, respectively. Finally, Section 6 covers the Faber--Krahn inequality for hinged external conditions and a double sided estimate for the clamped eigenvalue.

\begin{comment}
The positive number $\Lambda^H_{s,\infty}(\Om)$ (resp. $\Lambda^C_{s,\infty}(\Om)$) can be interpreted as the fractional generalization of $\Lambda^H_{\infty}(\Om)$ (resp. $\Lambda^C_{\infty}(\Om)$) in \cite{KP}. \\
\end{comment}

\section{Notation}\label{sec:not}
In the whole paper, $s\in(0,1)$ and $n\in\N$ are fixed numbers. \\
For $k\in\N\cup\{\infty\}$, we denote by $C_c^k(\R^n)$ and $C_0^k(\R^n)$ respectively the space of $C^k$-functions with compact support and the space of $C^k$-functions vanishing at infinity with all their derivatives up to order $k$. If $k=0$, we write simply $C_c(\R^n)$ and $C_0(\R^n)$. We denote by $\mathscr L^n$ the Lebesgue measure on $\R^n$ and for a measurable set $E\subseteq\R^n$ we write $|E|:=\mathscr{L}^n(E)$.\\
 Let $U\subseteq\R^n$ be an open set and $\alpha\in(0,1)$. We denote by $C^{\alpha}(U)$ the space of H\"older continuous functions of order $\alpha$ in $U$. By $C^{\alpha}_{\mathrm{loc}}(U)$ we refer to the space of locally H\"older continuous functions of order $\alpha$ in $U$, i.e. we say that $u\in C^{\alpha}_{\mathrm{loc}}(U)$ if for every $x\in U$ there exists a neighborhood $B$ of $x$ such that $u\in C^{\alpha}(B)$. More in general, if $\alpha\in\R^+\setminus\N$, we denote by $C^\alpha(U)$ (resp. $C_{\mathrm{loc}}^\alpha(U)$ ) the space of $\lfloor\alpha\rfloor$ times (resp. locally) differentiable functions with derivatives of order $\lfloor\alpha\rfloor$ being (resp. locally) H\"older continuous of order $\alpha-\lfloor\alpha\rfloor$.\\
 In the above, $(-\Delta)^su$ is the $s$-Laplacian of $u$, which is defined by
$$
(-\Delta)^s u(x):=c_{n,s}\,\int_{\R^n}\frac{u(x)-u(y)}{|x-y|^{n+2s}}dy\qquad\forall x\in\R^n,
$$
for $u$ smooth in $\R^n$. 
The previous integral is understood in the Cauchy principal value sense and $c_{n,s}:=\frac{s2^{2s\Gamma(\frac n2+s)}}{\pi^{n/2}\Gamma(1-s)}$, where $\Gamma$ denotes the Euler Gamma function. With this choice, the operator $(-\Delta)^s$ approaches the classical Laplacian as $s\to1^-$.\\
For $u\in L^\infty(\R^n)$, we interpret $(-\Delta)^su$ in the sense of distributions.\\
Suppose that $\gamma>0$ is such that $2s+\gamma$ is not an integer and $u\in C^{2s+\gamma}_{\mathrm{loc}}(\R^n)\cap L^1_s(\R^n)$. 
The space $L^1_s(\R^n)$ is the set of Lebesgue measurable functions $u:\R^n\to\R$ with 
$$
\int_{\R^n}\frac{|u(x)|}{1+|x|^{n+2s}}dx<\infty.
$$
Then, the $s$-Laplacian of $u$ is defined pointwise at every $x\in\R^n$ as
$$
(-\Delta)^s u(x):=c_{n,s}\,\mathrm{P.V.}\,\int_{\R^n}\frac{u(x)-u(y)}{|x-y|^{n+2s}}dy:=c_{n,s}\lim_{\eps\to0}\,\int_{\R^n\setminus B_\eps(x)}\frac{u(x)-u(y)}{|x-y|^{n+2s}}dy,
$$
where
$c_{n,s}:=\frac{s2^{2s\Gamma(\frac n2+s)}}{\pi^{n/2}\Gamma(1-s)}$ and $\Gamma$ denotes the Euler Gamma function.\\
Now assume that $u$ only belongs to $L^1_s(\R^n)$ and that $U$ is an open set of $\R^n$. We can define $(-\Delta)^su$ in the sense of distributions on $U$ as
$$
<(-\Delta)^su,\varphi>:=\int_{\R^n}u(-\Delta)^s\varphi\qquad\forall\varphi\in C^\infty_c(U).
$$
Notice that the integral in the right hand side is well defined since $(-\Delta)^s\varphi(x)\sim|x|^{-n-2s}$ away from $U$.\\
Moreover, if $p\in[1,\infty]$, we say that $(-\Delta)^su\in L^p(U)$ if there exists $f\in L^p(U)$ such that
$$
<(-\Delta)^su,\varphi>=\int_{U}f\varphi\qquad\forall\varphi\in C^\infty_c(U).
$$
%The abuse of notation is justified as follows: whenever $(-\Delta)^s u$ happens to be defined pointwise almost everywhere in $U$ (i.e. the limit defining the principal value exists at a.e. $x\in U$) and it is an $L^1_{\mathrm{loc}}(U)$-function, then, by integrating by parts, the distributional definition reduces to the pointwise one.
For $1\leq p<\infty$ and an open set $U\subseteq\R^n$, we define the fractional Sobolev space $W^{s,p}(U)$ as
$$
W^{s,p}(U)=\left\{u\in L^p(U): [u]_{W^{s,p}(U)}<\infty\right\},
$$
where $[\cdot]_{W^{s,p}(U)}$ stands for the Gagliardo seminorm, namely
$$
[u]_{W^{s,p}(U)}:=\left(\int_U dx\int_U\frac{|u(x)-u(y)|^p}{|x-y|^{n+sp}}dy\right)^{\frac1p}.
$$
Also, $W^{s,p}(U)$ is a Banach space, with the norm
$$
\|u\|_{W^{s,p}(U)}:=\left(\|u\|^p_{L^p(U)}+[u]^p_{W^{s,p}(U)}\right)^{1/p}.
$$
When $p=2$ we use the notation $H^s(U)$ to refer to the Hilbert space $W^{s,2}(U)$.\\
If $2s>1$, write $2s=1+\sigma$. Then, for $p\in[1,\infty)$, we define $$W^{2s,p}(U):=\{u\in W^{1,p}(U): Du\in W^{\sigma,p}(U)\},$$ equipped with the corresponding suitable norm $\|\cdot\|_{W^{2s,p}(U)}$.\\
Finally, for a bounded open set $\Omega\subset\R^n$ and $p\in (1,\infty)$, we denote
\begin{align*}
&W_0^{s,p}(\overline \Om):=\{u\in W^{s,p}(\R^n): u=0 \mbox{ in }\R^n\setminus\Om\},\\
&W_0^{s,p}(\Om):=\overline{C^\infty_c(\Om)}^{W^{s,p}(\Om)}.
\end{align*}
In particular, we write $H_0^s(\overline\Om):=W^{s,2}_0(\overline\Om)$.
If $\Omega$ has continuous boundary, then $W_0^{s,p}(\overline \Om)=W_0^{s,p}(\Om)$ (see \cite[Theorem 6]{FSV}).\\
Observe that for $u\in H_0^s(\overline\Om) $ we can define $(-\Delta)^s u$ in the weak sense as an element of the dual space $H^{-s}(\Om):=H^{s}_0(\Om)^*$, namely
$$
\langle(-\Delta)^su,v\rangle:=\frac{c_{n,s}}2\int_{\R^n\times\R^n}\frac{(u(x)-u(y))(v(x)-v(y))}{|x-y|^{n+2s}}dxdy\qquad\forall v\in H_0^s(\overline\Om),
$$
and we say that $(-\Delta)^su\geq0$ in the weak sense in $\Om$ if 
$$
\langle(-\Delta)^su,v\rangle\geq0\qquad\forall v\in H_0^s(\overline\Om):\, v\geq0.
$$

\section{Preliminaries}\label{prelim}

\subsection{On the regularity of competitors in $\mathscr L^{2s,\infty}_0(\Om)$ }\label{reg comp}
In this first preliminary section, we discuss some improvements in regularity in $\mathscr L^{2s,\infty}_0(\Om)$. In fact, we show that functions in this space are actually H\"older continuous and belong to $H^s_0(\Om)$.
\begin{proposition}
Set $$\mathbb{X}^s(\Omega):=\{u\in C^s(\R^n): u=0\,\mbox{ in }\,\R^n\setminus\Om, (-\Delta)^s u\in L^\infty(\Om)\}\subset \mathscr L^{2s,\infty}_0(\Om).$$
Then $\mathbb{X}^s(\Omega)\subset H^s_0(\overline\Om)$.
\end{proposition}
\begin{proof}
Since $u=0$ in $\R^n\setminus\Omega$ then
$$
[u]^2_{H^s(\R^n)}=\iint_{Q_\Omega}\frac{|u(x)-u(y)|^2}{|x-y|^{n+2s}}dydx,
$$
where $Q_\Omega:=\R^{2n}\setminus(\Omega^c\times\Omega^c)$. By using the symmetry of the integrand we have that
$$
[u]^2_{H^s(\R^n)}=[u]^2_{H^s(\Omega)}+2\int_\Omega\int_{\Omega^c}\frac{|u(x)|^2}{|x-y|^{n+2s}}dydx.
$$
Set $d(x):=d(x,\partial\Omega)$, for every $x\in\Omega$ we have that $\Omega^c\subset B^c_{d(x)}(x)$, therefore
$$
\int_{\Omega^c}\frac{dy}{|x-y|^{n+2s}}\le \int_{B^c_{d(x)}(x)}\frac{dy}{|x-y|^{n+2s}}=\frac{C'_{n,s}}{d^{2s}(x)}.
$$
Since $u\in C^s(\R^n)$, then $\frac{u}{d^s}\in L^\infty(\Omega)\subset L^2(\Omega).$\\
We have to prove $u\in H^s(\Omega)$. To this purpose, let us show that $u\in C^\alpha_{\rm loc}(\Om)$ for all $\alpha<2s$: Consider $(\eta_\varepsilon)_{\varepsilon>0}$ a sequence of mollifiers with compact support in $\Omega$, and set  $u_\varepsilon:=u\ast \eta_\varepsilon$, for every $\varepsilon>0$. Then $u_\varepsilon\in C^\infty_c(\R^n)$. By Proposition 2.3 in \cite{RS}, for fixed $\Om''\Subset\Om'\Subset\Om$ and every $\alpha\in (0,2s)$, we have
\begin{align}\label{holder est u eps}
\|u_\eps\|_{C^\alpha(\Om'')}\leq C(\|(-\Delta)^su_\eps\|_{L^\infty(\Om')}+\|u_\eps\|_{L^\infty(\R^n)}),
\end{align}
for some constant $C=C(n,s,\alpha)>0$. From the properties of the convolution, for every $\eps>0$, we have
$$
\|u_\eps\|_{L^\infty(\R^n)}\leq\|u\|_{L^\infty(\R^n)}\|\eta_\eps\|_{L^1(\R^n)}=\|u\|_{L^\infty(\R^n)}.
$$
Moreover,
for $\eps$ sufficiently small, we have
$$
\|(-\Delta)^su_\eps\|_{L^\infty(\Om')}\leq\|(-\Delta)^su\|_{L^\infty(\Om)}\|\eta_\eps\|_{L^1(\R^n)}=\|(-\Delta)^su\|_{L^\infty(\Om)}.
$$
Passing to the limit as $\eps\to0$ in \eqref{holder est u eps}, the lower semicontinuity of $\|\cdot\|_{C^\alpha(\Om'')}$ and the arbitrariness of $\Om''$ ensure that $u\in C^\alpha_{\rm loc}(\Om)$, for every $\alpha\in(0,2s)$.\\
We conclude that $u\in H^s(\Om)$ from the following estimate:
\begin{align*}
\iint_{\Omega\times\Om}\frac{|u(x)-u(y)|^2}{|x-y|^{n+2s}}dydx&\leq\iint_{\Omega\times\Om}\frac{|x-y|^{2\alpha}}{|x-y|^{n+2s}}dydx\\
&=\iint_{\Omega\times\Om}\frac{1}{|x-y|^{n+2s-2\alpha}}dydx<\infty.
\end{align*}
\end{proof}

Actually, the following stronger result holds.
\begin{proposition}\label{distr to weak}
Let $\Om\subset\R^n$ be a bounded Lipschitz domain satisfying a uniform exterior ball condition.
Then $$ \mathscr L^{2s,\infty}_0(\Om)\subset H^s_0(\Om).$$
\end{proposition}
\begin{proof}
Let $u\in \mathscr L^{2s,\infty}_0(\Om)$ and set $f=(-\Delta)^s u\in L^\infty(\Om)$. By Lax-Milgram's Theorem, there exists a unique $v\in H^s_0(\overline\Om)$ weak solution to
\begin{equation}\nonumber
\left\{
    \begin{aligned}
    (-\Delta)^s v&=f \qquad&&\mbox{in }\Om,\\
    v&=0 &&\mbox{in }\R^n\setminus\Om.
    \end{aligned}
    \right.
\end{equation}
By Proposition 1.1 in \cite{RS}, we have that $v\in C^s(\R^n)$. In particular, the function
$$w=u-v\in L^\infty(\R^n)$$
is a distributional solution to
\begin{equation}\nonumber
\left\{
    \begin{aligned}
    (-\Delta)^s w&=0 \qquad&&\mbox{in }\Om,\\
    w&=0 &&\mbox{in }\R^n\setminus\Om.
    \end{aligned}
    \right.
\end{equation}
By applying the Maximum Principle in \cite[Theorem 1.1]{GH}\footnote{Assumption (16) in \cite[Theorem 1.1]{GH} is fulfilled since $w\in L^\infty(\Om)$ and $$|\{x\in\Om: \mathrm{dist}(x,\R^n\setminus\Om)<\eps\}|\leq C\eps$$ for a fixed constant $C>0$, because $\Om$ has Lipschitz boundary.}, we necessarily have $w=0$ in $\Om$, ensuring $u=v\in H^s_0(\overline\Om)$. Finally, the latter space coincides with $H^s_0(\Om)$ since $\Om$ has Lipschitz boundary.
\end{proof}

\begin{corollary}\label{space equiv}
Let $\Om\subset\R^n$ be a bounded Lipschitz domain satisfying a uniform exterior ball condition. Then $$ \mathscr L^{2s,\infty}_0(\Om)=\mathbb{X}^s(\Omega)\subset H^s_0(\Om).$$
\end{corollary}
\begin{proof}
It is enough to observe that the inclusion $$\mathscr L^{2s,\infty}_0(\Om)\subset\mathbb{X}^s(\Omega)$$
holds true thanks to H\"older regularity estimates for the fractional Dirichlet problem (\cite[Proposition 1.1]{RS}).
\end{proof}

\subsection{Elliptic estimates}
In this second preliminary section, we prove some elliptic estimates for distributional solutions to the fractional Dirichlet problem with $L^1$ data. We start by recalling the following Calder\'on-Zygmund type result, contained in \cite[Corollary 1.7]{AFLY}, for $L^1$-weak solution.

\begin{definition}\label{L1 weak def}
Let $\Om\subset\R^n$ be a bounded open set with $C^2$-boundary and $h\in L^1(\Om)$. A function $u\in L^1(\Om)$ is an $L^1$-weak solution to
\begin{equation}\nonumber
\left\{
\begin{aligned}
(-\Delta)^s u&=h\qquad&&\mbox{in }\Om,\\
u&=0&&\mbox{in }\R^n\setminus\Om,
\end{aligned}
\right.
\end{equation}
if 
\begin{equation}\nonumber
\left\{
\begin{aligned}
\int_{\R^n}u(-\Delta)^s \varphi&=\int_\Om h\varphi\qquad&&\forall \varphi\in\mathbb X^s(\Om),\\
u&=0&&\mbox{in }\R^n\setminus\Om.
\end{aligned}
\right.
\end{equation}
\end{definition}
Existence and uniqueness for $L^1$-weak solutions are proved in \cite{CV}. Moreover, we notice
that in \cite{AFLY} $L^1$-weak solutions are simply called {\it{weak solutions}}.

\begin{proposition}[Fractional Calder\'on-Zygmund estimates]
\label{CZ thm}
Let $\Om\subset\R^n$ be a bounded domain with boundary of class $C^2$ and $s\in(0,1)$. Let $h\in L^m(\Om)$ and $u$ be the unique $L^1$-weak solution to 
\begin{equation}\nonumber
\left\{
\begin{aligned}
(-\Delta)^s u&=h\qquad&&\mbox{in }\Om,\\
u&=0&&\mbox{in }\R^n\setminus\Om.
\end{aligned}
\right.
\end{equation}
Then, we have:
\begin{itemize}
\item if $1\leq m< \frac ns$, then for all $1<p<m^*_s:=\frac{nm}{n-ms}$ there exists $C>0$ such that
$$
\|u\|_{W^{s,p}(\R^n)}\leq C\|h\|_{L^m(\Om)};
$$
\item if $m>\frac ns$, then for all $1<p<\infty$ there exists $C>0$ such that
$$
\|u\|_{W^{s,p}(\R^n)}\leq C\|h\|_{L^m(\Om)}.
$$
\end{itemize}
Here, $C > 0$ are positive constants depending only on $n, s, p, m$ and $\Om$.
\end{proposition}

\begin{remark}\label{rem CZ}
We observe that, since $m^*_s>m$, we can choose $p=m$ in Proposition \ref{CZ thm}, obtaining an estimate of the type $W^{s,p}(\R^n)-L^p(\Om)$ for all $p\in (1,\infty)$:
$$
\|u\|_{W^{s,p}(\R^n)}\leq C\|h\|_{L^p(\Om)},
$$
for a constant $C=C(n,s,p,\Om)>0$.
\end{remark}

Now, we introduce the notion of {\it{very weak solution}} to the Fractional Laplace problem. First, we define the space
$$
L^\infty_s(\R^n):=\{u\in L^\infty(\R^n): \|(1+|\cdot|^{n+2s})|u|\|_{L^\infty(\R^n)}<\infty\}.
$$
\begin{definition}\label{very weak def}
Let $\Om\subset\R^n$ be open, $g\in L^1_s(\R^n)$ and $h\in L^1(\Om)$. A function $u\in L^1_s(\R^n)$ is a very weak solution to
\begin{equation}\nonumber
\left\{
\begin{aligned}
(-\Delta)^s u&=h\qquad&&\mbox{in }\Om,\\
u&=g&&\mbox{in }\R^n\setminus\Om,
\end{aligned}
\right.
\end{equation}
if for every $\varphi\in L^\infty(\Om)$ with compact support in $\Om$, such that $(-\Delta)^s\varphi\in L^\infty_s(\R^n)$, it holds
\begin{equation}\nonumber
\left\{
\begin{aligned}
\int_{\R^n}u(-\Delta)^s \varphi&=\int_\Om h\varphi,\\
u&=g\qquad\mbox{in }\R^n\setminus\Om.
\end{aligned}
\right.
\end{equation}
\end{definition}

\begin{remark}\label{rem on class}
The test class in the above definition is contained in $\mathbb X^s(\Om)$, so Definition \ref{very weak def} is weaker (and then more general) than Definition \ref{L1 weak def} even in the case $g=0$.\\
Furthermore, the regularity of the tests $\varphi$ is actually $W^{2s,p}(\R^n)$ for all $p\in [2,\infty)$ (see \cite[Remark 3.5]{CM}). In particular, $(-\Delta)^s\varphi$ has a pointwise representative almost everywhere in $\R^n$.
\end{remark}

Now we recall a result contained in \cite{CCLP} about the regularity of very weak solutions to the homogeneous Fractional Dirichlet problem. We may state the result for any open set $\Om\subset\R^n$. 
\begin{proposition}[Regularity of very weak solutions with homogeneous data]\label{thm CCPL}
Let $\Om\subset\R^n$ be an open set and $s\in(0,1)$. Let $g\in L^1_s(\R^n)$ and let $u$ be a very weak solution to 
\begin{equation}\nonumber
\left\{
\begin{aligned}
(-\Delta)^s u&=0\qquad&&\mbox{in }\Om,\\
u&=g&&\mbox{in }\R^n\setminus\Om.
\end{aligned}
\right.
\end{equation}
Then, we have:
\begin{itemize}
\item[i)] $u\in H^s_\mathrm{loc}(\Om)$ and
$$
\|u\|_{H^s(B)}\leq C(B)\|u\|_{L^1_s(\R^n)}\qquad\forall B\Subset\Om;
$$
\item[ii)] $u$ is real analytic in $\Om$. Moreover, for any $\alpha\in \N^n_0$ and any $B'\Subset B\Subset\Om$, there holds
$$
\|D^\alpha u\|_{L^\infty(B')}\leq c^{|\alpha|}\alpha!C(B,B',n,s)\left(\|u\|_{L^\infty(B')}+\|u\|_{L^1_s(\R^n)}\right)
$$
for a universal constant $c>0$.
\end{itemize}
\end{proposition}

The main result of this section is the following.
\begin{proposition}[Sobolev regularity of very weak solutions with $L^1$ data]\label{Sob est very weak}
Let $\Om\subset\R^n$ be an open set and $s\in(0,1)$. Let $g\in L^1_s(\R^n)$ and $h\in L^1(\Om)$. Suppose that $u$ is a very weak solution to 
\begin{equation}\nonumber
\left\{
\begin{aligned}
(-\Delta)^s u&=h\qquad&&\mbox{in }\Om,\\
u&=g&&\mbox{in }\R^n\setminus\Om.
\end{aligned}
\right.
\end{equation}
Then for all $q\in [1,\frac{n}{n-s})$ and $\Om'\Subset\Om$ there holds
\begin{equation}\label{Sob reg est}
\|u\|_{W^{s,q}(\Om')}\leq C\left(\|u\|_{L^1_s(\R^n)}+\|h\|_{L^1(\Om)}\right)
\end{equation}
for a constant $C=C(n,s,q,\Om',\Om)>0$.
\end{proposition}
\begin{proof}
Let us fix an open subset $\Om''\Subset\Om$ with $C^2$-boundary and let $\psi$ be the unique $L^1$-weak solution to 
\begin{equation}\nonumber
\left\{
\begin{aligned}
(-\Delta)^s \psi&=h\qquad&&\mbox{in }\Om'',\\
\psi&=0&&\mbox{in }\R^n\setminus\Om''.
\end{aligned}
\right.
\end{equation}
From Proposition \ref{CZ thm}, for all $q\in[1,\frac{n}{n-s})$, we have
\begin{equation}\label{CZ for psi}
\|\psi\|_{W^{s,q}(\R^n)}\leq C\|h\|_{L^1(\Om)}
\end{equation}
for a constant $C=C(n,s,q,\Om'')>0$.\\
Now, set $f:=u-\psi$. Then $f$ is a very weak solution to 
\begin{equation}\nonumber
\left\{
\begin{aligned}
(-\Delta)^s f&=0\qquad&&\mbox{in }\Om'',\\
f&=u&&\mbox{in }\R^n\setminus\Om''.
\end{aligned}
\right.
\end{equation}
 From Proposition \ref{thm CCPL} ii), we have that $f$ is smooth in $\Om''$, so it is $s$-harmonic in $\Om''$ in the classical sense. In particular, it is a local weak solution to the Fractional Laplace problem in $\Om''$ (w.r.t. its own Dirichlet boundary data). So, from Theorems 1.4 and 3.2 in \cite{BLS}, for any open set $\Om'\Subset\Om''$, there exist $\sigma\in(0,1)$ and a constant $C=C(\sigma,n,s,\Om')$ such that 
$$
\|f\|_{C^{0,\sigma}(\Om')}\leq C\left(\|f\|_{L^2(\Om')}+\|f\|_{L^1_s(\R^n)}\right).
$$
Hence, by Proposition \ref{thm CCPL} i),
\begin{align*}
\|f\|_{C^{0,\sigma}(\Om')}\leq C\|f\|_{L^1_s(\R^n)}.
\end{align*}
Finally, using the above estimate and Proposition \ref{thm CCPL} ii), we have
$$
\|f\|_{W^{1,\infty}(\Om')}\leq C\left(\|f\|_{L^\infty(\Om')}+\|f\|_{L^1_s(\R^n)}\right)\leq C\|f\|_{L^1_s(\R^n)}
$$
for a constant $C=C(n,s,\Om',\Om'')>0$.\\
For all $q\in[1,\frac{n}{n-s})$, the triangle inequality gives
$$
\|u\|_{W^{s,q}(\Om')}-\|\psi\|_{W^{s,q}(\Om')}\leq C\left(\|u\|_{L^1_s(\R^n)}+\|\psi\|_{L^1(\Om)}\right).
$$
Estimate \eqref{Sob reg est} directly follows from the above one by applying \eqref{CZ for psi}.
\end{proof}

\section{Proof of Theorem \ref{thm1}}
In this section we prove Theorem \ref{thm1}. We divide the proof into two parts: first, we ensure the non triviality of $\Lambda^H_{s,\infty}(\Om)$; second, we show that every minimizer of the eigenvalue problem \eqref{eq:lambda} must be a multiple of the solution to the Dirichlet problem \eqref{talenti}. In particular, this automatically ensures the existence and uniqueness up to rescaling of a minimizer $u_\infty$, since \eqref{talenti} admits a unique weak solution $v$, which necessarily belongs to $\mathscr L^{2s,\infty}_0(\Om)$.

\subsection{Strict positivity of the minimum}\label{subsec strict}
Let us prove that $\Lambda^H_{s,\infty}(\Om)>0$.\\ For the sake of contradiction, suppose that $\Lambda^H_{s,\infty}(\Om)=0$. Then there exists a sequence $(u_k)_k\subset\mathscr{L}^{2s,\infty}_{0}(\Om)\setminus\{0\}$ with $$\frac{\|(-\Delta)^su_k\|_{L^\infty(\Om)}}{\|u_k\|_{L^\infty(\R^n)}}\to 0\qquad \mbox{as } k\to\infty.$$
Up to multiplying by a scaling factor, we can assume that $\|u_k\|_{L^\infty(\R^n)}=1$ for every $k\in\N$.
Thus, 
\begin{align}\label{fk}
(-\Delta)^su_k\to0\qquad\mbox{ uniformly in } \Om.
\end{align}
Set $f_k:=(-\Delta)^su_k\in L^\infty(\Om)$. Now, thanks to Proposition \ref{distr to weak}, we have that $u_k$ is a weak solution to
\begin{equation}\nonumber
\left\{
\begin{aligned}
(-\Delta)^su_k&=f_k\qquad&&\mbox{in }\Om,\\
u_k&=0 &&\mbox{in }\R^n\setminus\Om.
\end{aligned}
\right.
\end{equation}
 By Proposition 1 in \cite{RS}, we have
\begin{align*}
\|u_k\|_{C^{s}(\R^n)}&\leq C\|f_k\|_{L^\infty(\Om)}
\end{align*}
for a constant $C>0$ independent of $k$. But from \eqref{fk}, we have $\|f_k\|_{L^\infty(\Om)}\to 0$ as $k\to\infty$, so that
$$
u_k\to 0\qquad\mbox{uniformly in }\R^n,
$$
contradicting the fact that $\|u_k\|_{L^\infty(\R^n)}=1$ for all $k\in \N$.

\subsection{Characterization of minimizers}

\begin{proposition}
    Let $\Om\subset\R^n$ be a bounded Lipschitz domain satisfying a uniform exterior ball condition. Then $u\in\mathscr{L}^{2s,\infty}_0(\Omega)$ is a minimizer of $\Lambda^H_{s,\infty}(\Om)$ if and only if $u$ is a multiple of the torsion function, namely, a solution to \eqref{talenti}.
\end{proposition}

\begin{proof}
$``\Rightarrow"$
    Let $u\in\mathscr{L}^{2s,\infty}_0(\Omega)$ be a minimizer for $\Lambda^H_{s,\infty}(\Om)$. By Corollary \ref{space equiv}, we have $u\in H^s_0(\Om)\cap C^s(\overline\Om)$. With no loss of generality we assume that $\left\|(-\Delta)^s u\right\|_{L^\infty(\Omega)}=1$ and, by replacing $u$ with $-u$, that
    $$
    \|u\|_{L^\infty(\Om)}=\max_{\overline \Om}u.
    $$
    Now, let $v$ be a (classical) solution to \eqref{talenti}. Set $$g:=(-\Delta)^s u,\qquad w:=v-u\in H^s_0(\Om)\cap C^s(\overline\Om).$$ For sake of contradiction, suppose $g\not\equiv 1$. Thus, since $-1\leq g\leq 1$ a.e. in $\Om$, there exists a subset of $\Om$ of positive measure where $1-g>0$. We have that $w$ (weakly) solves
    \begin{equation*}
    \left\{
        \begin{aligned}
            (-\Delta)^s w&=1-g\quad&&\text{in }\Omega, \\
            w&=0\quad&&\text{in }\R^n\setminus\Omega.
        \end{aligned}
        \right.
    \end{equation*}
   Since $(-\Delta)^s w=1-g\geq 0$ a.e. in $\Om$, we have that $(-\Delta)^s w\geq 0$ also in the weak sense in $\Om$. Hence, by the Maximum Principle for weak solutions (see \cite[Lemma 6]{SV}), we deduce $w\geq0$ in $\Om$. Moreover, by the Strong Maximum Principle for weak solutions (see \cite[Corollary 4.2]{MN}), since $(-\Delta)^sw\not\equiv0$ in $\Om$, we must have $w>\inf_{\R^n}w=0$ in $\Om$. So that $u< v$ in $\Om$. Therefore, since $u=v=0$ in $\R^n\setminus\Om$ and $u,v\in C^s(\overline\Om)$, we must have
    $$
    \left\|u\right\|_{L^\infty(\Omega)}<\left\|v\right\|_{L^\infty(\Omega)}.
    $$
    We obtain
    $$
    \frac{\|(-\Delta)^sv\|_{L^\infty(\Om)}}{\|v\|_{L^\infty(\R^n)}}=\frac1{\|v\|_{L^\infty(\R^n)}}<\frac1{\|u\|_{L^\infty(\R^n)}}=\frac{\|(-\Delta)^su\|_{L^\infty(\Om)}}{\|u\|_{L^\infty(\R^n)}},
    $$
    which leads to a contradiction with the minimality of $u$.

    $``\Leftarrow"$

    Arguing similarly, let $u,v\in \mathscr{L}^{2s,\infty}_0(\Omega)\setminus\{0\}$, with $u$ solution of \eqref{talenti}, and set $g:=(-\Delta)^s v$. Assume $\left\|g\right\|_{L^\infty(\Omega)}=1$, otherwise we replace $v$ with $\frac{v}{\left\|(-\Delta)^s v\right\|_{L^\infty(\Omega)}}$. Then $|g|\le 1$ a.e. in $\Omega$. Setting $w:=u-v$, $z:=u+v$ we have that $w$ and $z$ solve
         \begin{equation*}
    \left\{
        \begin{aligned}
            (-\Delta)^s w&=1-g\geq0\quad&&\text{in }\Omega, \\
            w&=0\quad&&\text{in }\R^n\setminus\Omega.
        \end{aligned}
        \right.
    \end{equation*}

 \begin{equation*}
    \left\{
        \begin{aligned}
            (-\Delta)^s z&=1+g\geq0\quad&&\text{in }\Omega, \\
            z&=0\quad&&\text{in }\R^n\setminus\Omega.
        \end{aligned}
        \right.
    \end{equation*}

    The Maximum Principle yields $w,z\geq0$ in $\Omega$, therefore $|v|\leq u$ in $\Om$, so that
    $$
    \left\|v\right\|_{L^\infty(\Omega)}\le \left\|u\right\|_{L^\infty(\Omega)}.
    $$
    Hence
    $$
    \frac{\left\|(-\Delta)^s u\right\|_{L^\infty(\Omega)}}{\left\|u\right\|_{L^\infty(\Omega)}}=\frac{1}{\left\|u\right\|_{L^\infty(\Omega)}}\le \frac{1}{\left\|v\right\|_{L^\infty(\Omega)}}=\frac{\left\|(-\Delta)^s v\right\|_{L^\infty(\Omega)}}{\left\|v\right\|_{L^\infty(\Omega)}},
    $$
    which proves the minimality of $u$. If $u$ is a multiple of a solution to \eqref{talenti}, the claim follows by the scale invariance of $\Lambda^H_{s,\infty}(\Om)$.
    
\end{proof}

\section{Proof of Theorem \ref{thm2}}
We outline here the proof of Theorem \ref{thm2}. We split the argument into two subsections: in the first one, we show the existence of a minimizer $u^C_\infty$ via $L^p$-approximation; in the second one, we derive the PDE system \eqref{PDE clamped}. Notice that the non triviality of $\Lambda^C_{s,\infty}(\Om)$ directly follows from the one of $\Lambda^H_{s,\infty}(\Om)$, since clearly $$\Lambda^H_{s,\infty}(\Om)\leq\Lambda^C_{s,\infty}(\Om).$$

\subsection{Existence of a minimizer}\label{exist sec}
\begin{proof}
For simplicity, let us denote $\Lambda_\infty:=\Lambda^C_{s,\infty}(\Om)>0$. 
We also denote by $\lambda$ the normalized Lebesgue measure on $\Om$, in symbols $\lambda=\mathscr{L}^n/|\Om|$, which naturally induces the space
$$
L^{p}_\lambda(\Om)=\left\{u:\Om\to\R: \,\|u\|_{L^p_\lambda(\Om)}:=\left(\int_\Om |u|^p d\lambda\right)^\frac1p<\infty\right\}.
$$
We aim at showing existence of a minimizer for the supremal functional
$$
E_\infty(u):=\|(-\Delta)^su\|_{L^\infty(\R^n)}
$$
on the space $$\tilde{\mathscr{W}}^{2s,\infty}_0(\Om):=\{u\in\mathscr{W}^{2s,\infty}_0(\Om): \|u\|_{L^\infty(\R^n)}=1\}.$$
To this purpose, let us fix $p\geq2$ and $w\in C^2(\R^n)\cap C_0(\R^n)$, with $w>0$ and $\int_{\R^n}w\,dx=1$. Define the weighted potential space
$$
\mathscr{W}^{2s,p}_{w}(\R^n):=\left\{u\in L^{p}(\R^n):\int_{\R^n}\left|(-\Delta)^su\right|^pw\,dx<\infty\right\},
$$ 
equipped with the norm
$$
\|u\|_{\mathscr{W}^{2s,p}_w(\R^n)}:=\|u\|_{L^{p}(\R^n)}+\left(\int_{\R^n}\left|(-\Delta)^su\right|^pw\right)^{\frac1p}.
$$
In the above, $(-\Delta)^s u$ stands for the distributional fractional Laplacian of $u$.\\
We observe that $\mathscr{W}^{2s,p}_{w}(\R^n)$ is a reflexive Banach space, since it is isometrically isomorphic to a closed subspace of $L^{p}(\R^n)\times L^p(\R^n)$. We want to consider its closed subspace
$$
\mathscr{W}^{2s,p}_{w,0}(\Om):=\{u\in\mathscr{W}^{2s,p}_w(\R^n):u=0\quad \mbox{in }\R^n\setminus\Om\}.
$$
We notice that
\begin{align}\label{inclusion Sob}
\mathscr{W}^{2s,\infty}_0(\Om)\subseteq\mathscr{W}^{2s,p}_{w,0}(\Om)\subseteq W^{2s,p}(\R^n).
\end{align}
Indeed, the first inclusion is clear. Let us show the second one: for any $u\in \mathscr{W}^{2s,p}_{w,0}(\Om)$, there holds $(-\Delta)^s u\in L^p_{loc}(\R^n)$; moreover, since $u=0$ in $\R^n\setminus \Om$, we have that $(-\Delta)^s u$ is smooth outside $\Om$ and decays as $|x|^{-n-2s}$ at infinity, so that $(-\Delta)^s u\in L^p(\R^n)$. Thus $u$ belongs to the potential space $\mathscr L^{2s,p}(\R^n)$.  Since $p\geq2$, the latter is included in $W^{2s,p}(\R^n)$ (see \cite[Theorem 4.55]{ADV}). \\
Now we aim at minimizing the functional
\begin{align}\label{Ep}
E_p(u):=\left(\int_{\R^n}\left|(-\Delta)^su\right|^pw\,dx\right)^{\frac{1}{p}}
\end{align}
on $$\tilde{\mathscr{W}}^{2s,p}_{w,0}(\Om):=\{u\in{\mathscr{W}}^{2s,p}_{w,0}(\Om): \|u\|_{L_\lambda^p(\Om)}=1\}.$$
Let us prove that this is a weakly closed subspace of ${\mathscr{W}}^{2s,p}_{w,0}(\Om)$: let $(u_k)_k\subset\tilde{\mathscr{W}}^{2s,p}_{w,0}(\Om)$ be such that $u_k\rightharpoonup u$ weakly in ${\mathscr{W}}^{2s,p}_{w,0}(\Om)$ for some $u\in{\mathscr{W}}^{2s,p}_{w,0}(\Om)$. We want to show that $\|u\|_{L^p_\lambda(\Om)}=1$. To this purpose, set $f_k:=(-\Delta)^s u_k$. Then $f_k\in L^p_{loc}(\R^n)$ and, thanks to \eqref{inclusion Sob}, $u_k$ is a weak solution to 
     \begin{equation*}
    \left\{
        \begin{aligned}
            (-\Delta)^s u_k&=f_k\quad&&\text{in }\Omega', \\
            u_k&=0\quad&&\text{in }\R^n\setminus\Omega'.
        \end{aligned}
        \right.
    \end{equation*}
for any fixed domain $\Om'$ of class $C^2$ with $\Om\Subset\Om'\Subset\R^n$. By applying the Calderon-Zygmund estimates in Proposition \ref{CZ thm} (see also Remark \ref{rem CZ}), we can find a constant $C=C(\Om',s,n,p)$ such that
$$
\|u_k\|_{W^{s,p}(\R^n)}\leq C\|f_k\|_{L^p(\Om')}.
$$
On the other hand, we have
$$
\|f_k\|_{L^p(\Om')}\leq\left(\int_{\R^n}|(-\Delta)^su_k|^p\frac{w}{\inf_{\Om'}w}dx\right)^{\frac1p}\leq \frac{1}{\inf^{1/p}_{\Om'}w}\|u_k\|_{{\mathscr{W}}^{2s,p}_{w}(\Om)},
$$
so that
$$
\|u_k\|_{W^{s,p}(\R^n)}\leq C'\|u_k\|_{{\mathscr{W}}^{2s,p}_{w}(\Om)}
$$
for a constant $C'=C'(w,\Om',s,n,p)$. The weak compactness of the space $W^{s,p}(\R^n)$ ensures
$$
u_k\rightharpoonup u\quad\mbox{weakly in }W^{s,p}(\R^n).
$$
In particular, the above convergence takes place in $W^{s,p}(\Om)$, thus the compact Sobolev embeddings in Lipschitz domain allow us to conclude
$$
u_k\to u\quad\mbox{in }L^p(\Om).
$$
Since $\|u_k\|_{L^p_\lambda(\Om)}=1$ for all $k$, we must have $\|u\|_{L^p_\lambda(\Om)}=1$, which implies $u\in \tilde{\mathscr{W}}^{2s,p}_{w,0}(\Om)$.\\
Using direct methods, it is enough to show that $E_p$ is weakly lower semicontinuous and weakly coercive in ${\mathscr{W}}^{2s,p}_{w,0}(\Om)$. \\
First, we observe that $E_p$ is finite in $\tilde{\mathscr{W}}^{2s,p}_{w,0}(\Om)$. Indeed, let $\varphi\in C^\infty_c(\Om)\setminus\{0\}$ and set $u:=\varphi/\|\varphi\|_{L^p_\lambda(\Om)}$. Then $u\in \tilde{\mathscr{W}}^{2s,p}_{w,0}(\Om)$ and clearly $E_p(u)<\infty$.\\
We begin by showing the lower semicontinuity: let  $u_k\rightharpoonup u$ weakly in ${\mathscr{W}}^{2s,p}_{w,0}(\Om)$, then we have $E_p(u_k)^p\leq C$ for a constant $C$ independent of $k$. Therefore, $(-\Delta)^s u_k \,w^\frac1p\rightharpoonup(-\Delta)^s u\,w^\frac1p$ weakly in $L^p(\R^n)$. From the lower semicontinuity of the $\|\cdot\|_{L^p(\R^n)}$-norm, we conclude
$$
E_p(u)\leq\liminf_{k\to\infty} E_p(u_k).
$$
Concerning the coercivity, let us suppose that $E_p(u)\leq M$ for some $u\in{\mathscr{W}}^{2s,p}_{w,0}(\Om)$. Then, arguing as before, for any  $C^2$ domain $\Om'$ with $\Om\Subset\Om'\Subset\R^n$, the Calderon-Zygmund estimates imply
$$
\|u\|_{W^{s,p}(\R^n)}\leq C E_p(u)\leq CM
$$
for a constant $C=C(w,\Om',s,n,p)$.
In particular, we have $$\|u\|_{{\mathscr{W}}^{2s,p}_{w}(\Om)}\leq C(M+1).$$
Therefore, we have proved the existence of a minimizer $u_p\in \tilde{\mathscr{W}}^{2s,p}_{w,0}(\Om)$ for $E_p$ in $\tilde{\mathscr{W}}^{2s,p}_{w,0}(\Om)$.\\
Now fix $q\geq2$ and consider the sequence of (extended) functionals $(E_p)_p$ in $\mathscr{W}^{2s,q}_{w,0}(\Om)$ defined as
\begin{equation}\nonumber
E_p(u)=\left\{
\begin{aligned}
&\left(\int_{\R^n}|(-\Delta)^su|^pw\,dx\right)^{\frac{1}{p}}\!\!\!\!\!&,\quad&u\in \mathscr{W}^{2s,q}_{w,0}(\Om)\cap\mathscr{W}^{2s,p}_{w,0}(\Om),\\
&+\infty&,\quad&u\in \mathscr{W}^{2s,q}_{w,0}(\Om)\setminus \mathscr{W}^{2s,p}_{w,0}(\Om),
\qquad\end{aligned}
\right.
\begin{aligned}
\forall p\geq2.
\end{aligned}
\end{equation}
In a similar way, we can define $E_\infty$ in $\mathscr{W}^{2s,q}_{w,0}(\Om)$ as
\begin{equation}\nonumber
E_\infty(u)=\left\{
\begin{aligned}
&\|(-\Delta)^su\|_{L^\infty(\R^n)}\!\!\!\!\!&,\quad&u\in\mathscr{W}_{0}^{2s,\infty}(\Om),\\
&+\infty&,\quad&u\in \mathscr{W}^{2s,q}_{w,0}(\Om)\setminus\mathscr{W}_{0}^{2s,\infty}(\Om).
\end{aligned}
\right.
\end{equation}
  Using the H\"older inequality, it is not difficult to prove that $(E_p)_p$ is monotone increasing and
$$
\lim_{p\to\infty}E_p(u)=E_\infty(u)\qquad\forall u\in\mathscr{W}^{2s,q}_{w,0}(\Om).
$$
Let us show that $(E_p)_p$ is equi-coercive in $\mathscr{W}^{2s,q}_{w,0}(\Om)$ with respect to $p$. Assume $E_p(u)\leq M$ for all $p\geq q$. By fixing a domain $\Om'\Supset\Om$ with $C^2$-boundary, we can apply the fractional Calderon-Zygmund estimates as before. Together with the H\"older inequality, this yields
$$
\|u\|_{W^{s,q}(\R^n)}\leq C E_q(u)\leq CE_p(u)\leq CM
$$
for a constant $C=C(\Om',s,n,q,w)$. Thus,
$$\|u\|_{{\mathscr{W}}^{2s,q}_{w}(\Om)}\leq C(M+1),$$
giving the equi-coercivity of $(E_p)_p$. \\
So, $(E_p)_p$ is a monotone increasing sequence of weakly lower semicontinuous functionals
 on $\mathscr{W}^{2s,q}_{w,0}(\Om)$. Therefore, its $\Gamma$-limit is given by (see \cite[Remark 2.12]{Br})
$$
\Gamma-\lim_p E_p=\lim_pE_p=E_\infty.
$$
Moreover, by the Fundamental Theorem of Gamma Convergence \cite[Thm. 2.10]{Br}, if $(u_p)_p\subset\tilde{\mathscr{W}}^{2s,p}_{w,0}(\Om)$ is a sequence of minimizers of $E_p$ in $\tilde{\mathscr{W}}^{2s,p}_{w,0}(\Om)$, up to passing to a subsequence, we have $u_p\rightharpoonup u_\infty$ weakly $\mathscr{W}^{2s,q}_{w,0}(\Om)$ where $u_\infty$ is a minimizer of $E_\infty$ on $\mathscr{W}_{w,0}^{2s,q}(\Om)$. By the definition of $E_\infty$ on $\mathscr{W}^{2s,q}_{w,0}(\Om)$, we necessarily have $u_\infty\in\mathscr{W}_{0}^{2s,\infty}(\Om)$. Moreover, there holds $\|u_\infty\|_{L^\infty(\Om)}=1$. Indeed
using the H\"older inequality and the compact Sobolev embeddings, we have
$$
1=\|u_p\|_{L^p_\lambda(\Om)}\leq\|u_p\|_{L^\infty(\Om)}\to\|u_\infty\|_{L^\infty(\Om)}\qquad\mbox{as }p\to\infty;
$$
on the other hand, 
$$
\|u_\infty\|_{L^q_\lambda(\Om)}=\lim_p\|u_p\|_{L^q_\lambda(\Om)}\leq\limsup_p\|u_p\|_{L^p_\lambda(\Om)}=1.
$$
Thus, $u_\infty\in\tilde{\mathscr{W}}^{2s,\infty}_0(\Om)$.
Furthermore, Theorem 2.10 in \cite{Br} also implies
\begin{align}\label{limit ep}
E_p(u_p)\to E_\infty({u_\infty})=\Lambda_{\infty}\quad\mbox{as }p\to\infty.
\end{align}
\end{proof}

\subsection{PDE derivation}

\begin{proof}
In the previous section, for any $p\geq 2$, we proved the existence of a minimizer $u_p$ to $E_p$ in $\mathscr{W}^{2s,p}_{w,0}(\Om)$ subjected to the constraint
$$
G_p(u):= \frac 1p\left(\int_{\Omega}|u|^p\,d\lambda-1\right)=0.
$$
Let us denote $\Lambda_p:=E_p(u_p)$. 
Then, from the Lagrange multiplier Theorem, there exists $\sigma\in\R$ such that
\begin{align}\label{EL up}
D E_p(u_p)[\varphi]=\sigma DG_p(u_p)[\varphi]\qquad\forall \varphi\in\mathscr{W}^{2s,p}_{w,0}(\Om),
\end{align}
where
\begin{align*}
D E_p(u_p)[\varphi]&:=\frac{d}{dt}\Big|_{t=0} E_p(u_p+t\varphi)\\
&=\Lambda_p^{1-p}\int_{\R^n}|(-\Delta)^s u_p|^{p-2}(-\Delta)^s u_p\,w\,(-\Delta)^s \varphi\, dx
\end{align*}
and
\begin{align*}
D G_p(u_p)[\varphi]&:=\frac{d}{dt}\Big|_{t=0} G_p(u_p+t\varphi)=\int_{\Om}| u_p|^{p-2} u_p \varphi\, d\lambda.
\end{align*}
We observe that $\sigma=\Lambda_p$. Indeed, by testing \eqref{EL up} with $u_p$, we get
\begin{align*}
\sigma&=\sigma\int_{\Om}|u_p|^pd\lambda=\Lambda_p^{1-p}\int_{\R^n}|(-\Delta)^s u_p|^pw\,dx=\Lambda_p^{1-p}\Lambda_p^{p}=\Lambda_p.
\end{align*}
Setting 
\begin{align}\label{fp}
f_p:=\Lambda_p^{-p}|(-\Delta)^s u_p|^{p-2}(-\Delta)^s u_p\,w
\end{align}
and
$$
g_p:=\frac{|u_p|^{p-2}u_p}{|\Om|},
$$
the Euler-Lagrange equation in \eqref{EL up} reads as
\begin{align}\label{EL up short}
(-\Delta)^s f_p= g_p \qquad \mbox{in }\mathscr W^{2s,p}_{w,0}(\Om)^*.
\end{align}
Notice that $f_p$ is a well defined measurable function on $\R^n$, at least for sufficiently large $p$, i.e. $p\geq p_0$ for some $p_0\geq2$, since $\Lambda_p\to\Lambda_\infty>0$.\\
Now, using the H\"older inequality, we have
\begin{equation}\label{L1 bounds}
\left\{
\begin{aligned}
    &\|f_p\|_{L^1(\R^n)}\leq \Lambda_p^{-1}\\
    &\|g_p\|_{L^1(\R^n)}\leq1
\end{aligned}
\right.
\qquad\forall p\geq p_0.
\end{equation}
Indeed,
\begin{align*}
\|f_p\|_{L^1(\R^n)}&=\Lambda_p^{-p}\int_{\R^n}|(-\Delta)^su_p|^{p-1}w\,dx\\
&\leq\Lambda_p^{-p}\left(\int_{\R^n}|(-\Delta)^s u_p|^pw\right)^\frac{p-1}{p}\left(\int_{\R^n}w\right)^\frac1p=\Lambda_p^{-p}\Lambda_p^{p-1}=\Lambda_p^{-1},
\end{align*}
and similarly, recalling $\|u_p\|_{L^p_\lambda(\Om)}=1$,
\begin{align*}
\|g_p\|_{L^1(\R^n)}&=\frac1{|\Om|}\int_{\Om}|u_p|^{p-1} dx\\
&=\int_\Om|u_p|^{p-1}d\lambda\leq\left(\int_\Om|u_p|^pd\lambda\right)^\frac{p-1}{p}\left(\int_{\Om}d\lambda\right)^\frac1p=1.
\end{align*}
We observe that any test function in Definition \ref{very weak def} clearly belongs to $\mathscr W^{2s,p}_{w,0}(\Om)$. Thus, from \eqref{EL up short}, $f_p$ can be regarded as a very weak solution to the Dirichlet problem with $L^1$ data $g_p$, with respect to its own external datum in $\R^n\setminus\Om$. By Proposition \ref{Sob est very weak}, for all $q\in [1,\frac{n}{n-s})$ and $\Om'\Subset\Om$ there holds
\begin{equation}
\|f_p\|_{W^{s,q}(\Om')}\leq C\left(\|f_p\|_{L^1_s(\R^n)}+\|g_p\|_{L^1(\Om)}\right)
\end{equation}
for a constant $C=C(n,s,q,\Om',\Om)>0$.
Since $\Lambda_p\to\Lambda_\infty$ as $p\to\infty$, we have that both $(f_p)_p$ and $(g_p)_p$ are bounded in $L^1(\R^n)$, so that $(f_p)_p$ is bounded also in $W^{s,q}_{loc}(\Om)$ for all $q\in [1,\frac{n}{n-s})$. From the Sobolev compact embeddings (see \cite[Thm. 7.1]{DPV}), up to passing to a subsequence, we have 
\begin{align}\label{fp conv}
f_p\to f_\infty\qquad\mbox{ in } L^q_{loc}(\Om)
\end{align}
for a function $f_\infty\in W^{s,q}_{loc}(\Om)\cap L^1(\Om)$.\\
Moreover, the $L^1$-bounds in \eqref{L1 bounds} ensure the existence of measures $\mu,\nu\in\mathcal{M}(\R^n)$ with $|\mu|(\R^n)\leq\Lambda_\infty^{-1}$ and $|\nu|(\R^n)=|\nu|(\overline\Om)\leq1$ such that (up to passing to further subsequences)
\begin{equation}
\left\{
\begin{aligned}
    &f_p\stackrel{*}{\rightharpoonup}\mu\\
    &g_p\stackrel{*}{\rightharpoonup}\nu
\end{aligned}
\right.
\qquad\mbox{in }\mathcal{M}(\R^n)
    \end{equation}
In particular, from \eqref{fp conv}, we necessarily have $\mu\mres\Om=f_\infty\mathscr{L}^n\mres\Om$. 
Furthermore, since $u_p=0$ in $\R^n\setminus\Om$, the support of $\nu$ is clearly contained in $\overline\Om$. Actually, we can easily show that 
$$
\mathrm{supp}(\nu)\subseteq\mathrm{arg max}(|u_\infty|).
$$
Indeed, if $u_\infty(x)\in (-1,1)$, then from the uniform convergence of $u_p$ to $u_\infty$, there exists $\eps>0$, independent of $p$, such that $|u_p(x)|<1-\eps$ for all $p$. So that 
$$
|g_p(x)|=\frac{|u_p(x)|^{p-1}}{|\Om|}\to 0 \quad\mbox{as }p\to\infty,
$$
implying $x\notin\mathrm{supp}(\nu)$.
In particular, we have that supp$(\nu)\Subset\Om$.\\
Passing to the limit in \eqref{EL up short}, we obtain
\begin{align}\label{meas eq}
    (-\Delta)^s \mu=\nu\quad\mbox{in }\mathcal D'(\Om).
\end{align}
Let us show that $\mu,\nu,f_\infty\not\equiv0$. \\
Using the uniform convergence of $u_p$ to $u_\infty$, we have
$$
1=\int_\Om g_pu_pdx\to\int_{\overline \Om}u_\infty d\nu=\int_\Om u_\infty d\nu,
$$
so that $\nu\neq0$ and, from \eqref{meas eq}, also $\mu\neq0$.
It remains to show that $f_\infty\not\equiv0$. \\
For any $\eps\in \left(0,\frac{|\nu|(\Om)}{1+\Lambda_\infty^{-1}}\right)$, we can choose a function $\varphi\in C^0_c(\Om)$ with 
$$
\mathrm{supp}(\nu)\subseteq\mathrm{supp}(\varphi), \qquad\int_\Om\varphi \,d\nu>|\nu|(\Om)-\eps.
$$
By arguing by approximation in the uniform convergence in $\Om$, we can assume that $\varphi\in C^{2s+\gamma}_c(\Om)$ for some $\gamma>0$. Now, using \cite[Corollary 1]{CM}, we can modify $\varphi$ in a ball $B_r(x_0)\Subset\Om\setminus\mathrm{supp}(\nu)$ and construct a function $\tilde\varphi\in C^{2s+\gamma}_c(\Om)$ with 
\begin{equation}\nonumber
\left\{
\begin{aligned}
    &|(-\Delta)^s \tilde\varphi|\leq\eps&\quad  \mbox{in }\R^n\setminus B_r(x_0)\\
    &\tilde\varphi=\varphi& \mbox{in }\R^n\setminus B_r(x_0).
\end{aligned}
\right.
\end{equation}
Observe that $\tilde\varphi=\varphi$ on $\mathrm{supp}(\nu)$, so
$$
\int_\Om\tilde\varphi \,d\nu=\int_{\mathrm{supp}(\nu)}\varphi \,d\nu>|\nu|(\Om)-\eps.
$$
An approximation argument allows us to test \eqref{meas eq} with $\tilde\varphi$ and obtain
\begin{align*}
    |\nu|(\Om)-\eps&<\int_\Om\tilde\varphi\,d\nu\\
    &=\int_{\R^n} (-\Delta)^s\tilde\varphi\,d\mu\\
    &=\int_{\Om} f_\infty(-\Delta)^s\tilde\varphi\,dx+\int_{\R^n\setminus\Om} (-\Delta)^s\tilde\varphi\,d\mu\\
    &\leq\int_{\Om} f_\infty(-\Delta)^s\tilde\varphi\,dx+\eps|\mu|(\Om)\\
    &\leq\int_{\Om} f_\infty(-\Delta)^s\tilde\varphi\,dx+\eps\Lambda_\infty^{-1}.
\end{align*}
Therefore,
$$
\int_{\Om} f_\infty(-\Delta)^s\tilde\varphi\,dx>|\nu|(\Om)-(1+\Lambda_\infty^{-1})\eps>0,
$$
ensuring $f_\infty\not\equiv0$ in $\Om$.\\
We are left with proving \eqref{PDE clamped}. Recall that, for all $q\in [2,\infty)$, $u_p\rightharpoonup u_\infty$ weakly in $\mathscr{W}^{2s,q}_{w,0}(\Om)$. Thus, up to passing to a further subsequence, we have 
\begin{align}\label{weak conv Lapl}
(-\Delta)^su_p\rightharpoonup (-\Delta)^su_\infty\quad\mbox{ weakly in }L^q_{\mathrm{loc}}(\R^n).
\end{align}
In particular, the convergence in \eqref{weak conv Lapl} holds (weakly) in $L^q(\Om)$. Recalling \eqref{fp}, we can write
\begin{equation}
\label{eq:pdeapproxp}
|f_p|^\frac{1}{p-1}\mathrm{sgn}(f_p)=\Lambda_p^{-\frac{p}{p-1}}w^\frac{1}{p-1}(-\Delta)^su_p\qquad\mbox{in }\Om.
\end{equation}
If we restrict the above equation to a (non-empty) compact set $K\subset\Om\setminus\{f_\infty=0\}$ and pass to the limit as $p\to\infty$ in \eqref{eq:pdeapproxp}, recalling \eqref{fp conv}, we infer
$$
\mathrm{sgn}(f_\infty)=\Lambda_\infty^{-1}(-\Delta)^su_\infty\qquad\mbox{in }K,
$$
which provides \eqref{PDE clamped}.\\
Finally, one can prove \eqref{EL out Om} and that supp$(\mu)\Subset\R^n$ reasoning as in \cite[Sec. 4.3]{CM}.
\end{proof}

\section{Faber-Krahn inequality}

In this final section, we prove Theorem \ref{FK hinged}, concerning the Faber-Krahn inequality in the hinged case and its rigidity. Lastly, we exhibit a double sided estimate of the first eigenvalue for the clamped problem.

\begin{proof}
We want to perform the Talenti symmetrization argument proved in \cite{FV}. By Theorem \ref{thm1} we argue on the solutions to
 \begin{equation*}
    \left\{
        \begin{aligned}
            (-\Delta)^s u&=1\quad&&\text{in }\Omega, \\
            u&=0\quad&&\text{in }\R^n\setminus\Omega.
        \end{aligned}
        \right.
    \end{equation*}
and
 \begin{equation*}
    \left\{
        \begin{aligned}
            (-\Delta)^s v&=1\quad&&\text{in }B_R, \\
            v&=0\quad&&\text{in }\R^n\setminus B_R.
        \end{aligned}
        \right.
    \end{equation*}
Now we consider the Schwarz rearrangement $u^\sharp: B_R\to\R$ given by
$$
u^\sharp(x)=\sup\{t\geq0: |\{ u>t\}|>\omega_n|x|^n\},
$$
where $\omega_n$ is the measure of the unit ball in $\R^n$.
Then, the mass comparison (Theorem 3.1 in \cite{FV})
$$
\int_{B_r} u^\sharp (x)dx\le \int_{B_r} v(x)dx,
$$
for every $r>0$, induces
$$
\left\|u\right\|_{L^\infty(\Omega)}=\left\|u^\sharp\right\|_{L^\infty(B_R)}\le \left\|v\right\|_{L^\infty(B_R)},
$$
which implies that $\Lambda^H_{s,\infty}(\Omega)\ge \Lambda^H_{s,\infty}(B_R)$.
\end{proof}

\begin{corollary}\label{cor:u sulla palla}
Let $B_R$ be the ball centered at the origin of radius $R$. Then, every minimizer is a multiple of the fractional torsion function
$$
u_\infty^H(x):=K_{n,s}(R^2-|x|^2)^s_+,
$$
so that
$$
\Lambda^H_{s,\infty}(B_R)=\frac{K^{-1}_{n,s}}{R^{2s}},
$$
where $K_{n,s}=\frac{\Gamma\left(\frac n2\right)}{2^{2s}\Gamma(1+s)\Gamma\left(\frac n2+s\right)}$.
\end{corollary}

\medskip
\noindent\textit{Rigidity.}
We now prove the rigidity of the inequality \eqref{eq:FKH}. To do this, we firstly state and prove the following preliminary lemmas

\begin{lemma}
\label{lem:heat_symmetrization}
Let $D\subset\mathbb R^n$ be a bounded domain
and let $B:=B_R(x_0)$ be the ball centered at $x_0$ with $R>0$ such that $|D|=|B|$.\\
Denote by $p_D(t,x,y)$ and $p_{B}(t,x,y)$
the Dirichlet heat kernels associated with the
fractional Laplacian $(-\Delta)^s$ in $D$ and
$B$, respectively.\\
Then, for every $t>0$ and $x\in D$,
\begin{equation}\label{eq:heat_symmetrization}
\int_D p_D(t,x,y)\,dy
\leq
\int_{B}p_{B}(t,x_0,y)\,dy.
\end{equation}
Equivalently, denoting by
$$
h_D(t,x):=\int_D p_D(t,x,y)\,dy
$$
the corresponding heat distribution, one has
$$
h_D(t,x)\leq h_{B}(t,x_0)
\qquad\text{for every }t>0,\ x\in D.
$$
\end{lemma}

\begin{proof}
Let $(X_t)_{t\geq0}$ be the rotationally
symmetric $2s$-stable L\'evy process with
infinitesimal generator $-(-\Delta)^s$.
For an open set $U\subset\mathbb R^n$, set
$$
\tau_U:=\inf\{t>0:X_t\notin U\}.
$$
Since the process is rotationally symmetric,
its symmetrization coincides with the original
process. Thus, by the symmetrization inequality
in \cite[Theorem 1.4]{BMH10}, applied with
$f\equiv1$ and $V\equiv0$, we obtain
$$
\mathbb P_x(\tau_D>t)
\leq
\mathbb P_{x_0}(\tau_{B}>t)
\qquad\text{for every }t>0.
$$
Recalling that
$$
\mathbb P_x(\tau_U>t)
=\int_U p_U(t,x,y)\,dy,
$$
the desired inequality immediately follows.
\end{proof}

\begin{lemma}
\label{lem:smalltimeasymp}
Let $D$ be a bounded domain with $|\partial D|=0$ and let $h_D$ be the associated heat distribution, namely the solution to
\begin{equation}
\begin{cases}
\partial_t h_D+(-\Delta)^s h_D=0
& \text{in }(0,\infty)\times D,\\
h_D=0
& \text{in }(0,\infty)\times D^c,\\
h_D(0,\cdot)=1
& \text{in }D.
\end{cases}
\end{equation}
Then, for every $x\in D$,
\begin{equation}
\label{eq:lemmaheatkernel}
\lim_{t\to 0^+}\frac{1-h_D(t,x)}{t}=\kappa_D(x),
\end{equation}
where
$$
\kappa_D(x)=c_{n,s}\int_{D^c}
\frac{dy}{|x-y|^{n+2s}}.
$$
\end{lemma}
\begin{proof}
For every $t>0$ and $u_0\in L^1_{\rm loc}(\R^n)$, we define $$T_D(t)u_0:=\int_D p_D(t,\cdot,y)u_0(y)dy,$$ the fractional Dirichlet heat semigroup, so that $h_D(t,\cdot)=T_D(t)\chi_D$ in $D$.
For $\varphi\in C_c^\infty(\mathbb R^n)$, set
$g_\varphi=(-\Delta)^s\varphi|_D$ and define
\begin{equation}
\label{eq:duhamel}
  v_\varphi(t,\cdot)
  =\varphi-\int_0^t T_D(\tau)g_\varphi\,d\tau,
\end{equation}
where $T_D(\tau)g_\varphi$ is extended by zero to $D^c$.
By Duhamel's formula, $v_\varphi-\varphi$ solves the fractional heat equation
with source $-g_\varphi$, zero initial datum, and zero exterior
datum. Consequently, $v_\varphi$ solves
\[
  \partial_t v_\varphi+(-\Delta)^s v_\varphi=0
  \quad\text{in }(0,\infty)\times D,
\]
with initial datum $\varphi|_D$ and exterior datum $\varphi|_{D^c}$.\\
Fix $x\in D$ and $r>0$ such that $B_{2r}(x)\subset D$, and let $\varphi^-_j\in C^\infty_c(D)$, $\varphi^+_j\in C^\infty_c(\R^n)$ such that
$$
0\le \phi_j^-\le \chi_D\le \phi_j^+\le 1
\quad\text{in }\mathbb R^n,
\qquad
\phi_j^-=\phi_j^+=1\quad\text{on }B_r(x),
$$
and $\phi_j^\pm\to\chi_D$ pointwise in $\R^n\setminus\partial D$ as $j\to\infty$.
Comparison principle for solutions to parabolic equation gives
\begin{equation}
\label{eq:carabinieri}
  v_{\varphi^-_j}(t,x)\le h_D(t,x)\le v_{\varphi^+_j}(t,x).
\end{equation}
For every $\varphi_j^\pm$, the function
$g_{\varphi^\pm_j}$ is bounded and continuous. By combining \eqref{eq:duhamel} and the fact that $T_D(\cdot)$ is a Feller semigroup (see \cite[Lemma~2.2]{BK24}), it holds that
\[
  \lim_{t\to 0^+}\frac{1-v_{\varphi^\pm_j}(t,x)}{t}
  =\lim_{t\to 0^+}\frac1t\int_0^t T_D(\tau)g_{\varphi^\pm_j}(x)\,d\tau
  =g_{\varphi^\pm_j}(x)
  =(-\Delta)^s\varphi^\pm_j(x).
\]
Since every $\varphi_j^\pm$ equals one on $B_r(x)$
and $|\partial D|=0$, dominated convergence yields
$$
  (-\Delta)^s\varphi_j^\pm(x)
  =c_{n,s}\int_{\mathbb R^n}
      \frac{1-\varphi_j^\pm(y)}{|x-y|^{n+2s}}\,dy
  \longrightarrow \kappa_D(x)\qquad\mbox{as }j\to\infty.
$$
Thus, from \eqref{eq:carabinieri}, we conclude that
$$
 \lim_{t\to 0^+}\frac{1-h_D(t,x)}{t}=\lim_{j\to\infty}\lim_{t\to 0^+}\frac{1-v_{\varphi^\pm_j}(t,x)}{t}=\kappa_D(x).
$$
\end{proof}

\begin{proof}[\textit{Rigidity of the Faber-Krahn inequality for the first hinged eigenvalue}]
Let $u_\Omega$ be the eigenfunction relative to $\Lambda_{s,\infty}^H(\Omega)$ and let $x_0\in\Omega$ be a maximum point of
$u_\Omega$, so that
\[
u_\Omega(x_0)
=\|u_\Omega\|_{L^\infty(\Omega)}.
\]
Set $B:=B_R(x_0)$, where $R>0$ is chosen
such that $|B|=|\Omega|$.\\
Now, let us suppose that $\Omega$ does not coincide with
$B$ up to a set of measure zero. Since
$|\Omega|=|B|$, we necessarily have
\[
|B\setminus\Omega|
=
|\Omega\setminus B|>0.
\]
On the other hand, we observe that
\begin{align*}
&\kappa_\Omega(x_0)-\kappa_B(x_0)\\
&\quad =
c_{n,s}\left(
\int_{B\setminus\Omega}
\frac{dy}{|x_0-y|^{n+2s}}
-
\int_{\Omega\setminus B}
\frac{dy}{|x_0-y|^{n+2s}}
\right).
\end{align*}
Since
\[
|x_0-y|<R
\quad\text{for }y\in B\setminus\Omega,
\]
whereas
\[
|x_0-y|>R
\quad\text{for }y\in\Omega\setminus B,
\]
the decreasing monotonicity of
$r\mapsto r^{-n-2s}$, together with
$|B\setminus\Omega|=|\Omega\setminus B|>0$,
implies
\[
\kappa_\Omega(x_0)>\kappa_B(x_0).
\]
Combining this strict inequality with Lemma \ref{lem:smalltimeasymp}, we infer that
there exists $t_0>0$ such that
\begin{equation}
\label{eq:estimateheatdist1}
h_\Omega(t,x_0)<h_B(t,x_0)
\qquad\text{for every }t\in(0,t_0).
\end{equation}
On the other hand, by Lemma \ref{lem:heat_symmetrization}, we have
\begin{equation}
\label{eq:estimateheatdist2}
h_\Omega(t,x)\leq h_B(t,x_0)
\qquad\text{for every }x\in\Omega,\, t>0.
\end{equation}
Therefore, recalling the representation of
the fractional torsion function in terms of its heat
distribution, and putting together \eqref{eq:estimateheatdist1} and \eqref{eq:estimateheatdist2} we obtain
\begin{align*}
&u_B(x_0)-u_\Omega(x_0)\\
&=\int_0^\infty
\bigl(h_B(t,x_0)-h_\Omega(t,x_0)\bigr)\,dt\\
&=\int_0^{t_0}
\bigl(h_B(t,x_0)-h_\Omega(t,x_0)\bigr)\,dt\\
&\quad+\int_{t_0}^\infty
\bigl(h_B(t,x_0)-h_\Omega(t,x_0)\bigr)\,dt\\
&>0,
\end{align*}
since the first integral is strictly positive,
whereas the second is nonnegative.\\
Consequently, since by Corollary \ref{cor:u sulla palla} the expression of $u_B$ is explicit, we have
$$
\|u_\Omega\|_{L^\infty(\Omega)}=u_\Omega(x_0)
<
u_B(x_0)=\|u_B\|_{L^\infty(B)},
$$
which implies
$$
\Lambda^H_{s,\infty}(\Omega)
>
\Lambda^H_{s,\infty}(B).
$$

Thus, equality in the Faber--Krahn inequality
forces $\Omega$ to coincide with a ball
up to a $\mathcal{L}^n$-negligible set.
Conversely, equality clearly holds when
$\Omega$ is a ball, by translation invariance.
This completes the proof.
\end{proof}

We conclude by proving the announced Proposition \ref{prop:estimatesforclamp}. It remains open whether the ball is the optimal set even for the clamped eigenvalue.

\begin{proof}
As already observed, since $\mathscr L^{2s,\infty}_0(\Om)\supset\mathscr W^{2s,\infty}_0(\Om)$ for every admissible set $\Omega$, it holds that
$$
\Lambda^C_{s,\infty}(\Om)\ge \Lambda^H_{s,\infty}(\Om),
$$
and, in particular,
$$
\Lambda^C_{s,\infty}(B_R)\ge \Lambda^H_{s,\infty}(B_R)=\frac{K^{-1}_{n,s}}{R^{2s}},
$$
where $K_{n,s}$ is the same as in Corollary \ref{cor:u sulla palla}. For the upper estimate, let us consider the function $w(x):=(R^2-|x|^2)^{s+1}_+$. Since $w\in C^{1,s}_c(\R^n)$, to prove that $w\in \mathscr W^{2s,\infty}_0(B_R)$ we are left to prove that $(-\Delta)^s w\in L^\infty(\R^n)$. By \cite[Lemma 7 and Table 3]{Dyda}, for $x\in \overline{B}_R$, we have  
\begin{align}\label{Dyda formula}
(-\Delta)^sw(x)=(1+s)K^{-1}_{n,s}\left(R^2-\frac{n+2s}{n}|x|^2\right).
\end{align}
So that
$$
|(-\Delta)^s w|\le (1+s)K^{-1}_{n,s}R^2\max\left\{1,\frac{2s}{n}\right\}\qquad\mbox{in }\overline B_R.
$$
On the other hand, for $x\in\R^n\setminus\overline{B}_R$ we have
$$
|(-\Delta)^s w(x)|=c_{n,s}\int_{B_R}\frac{w(y)}{|x-y|^{n+2s}}dy.
$$
Since $(-\Delta)^s w$ is radially symmetric in $\R^n$, from the above identity we deduce that $|(-\Delta)^s w|$ is radially decreasing in $\R^n\setminus\overline{B}_R$. In particular, 
since $(-\Delta)^s w\in C^{0,1-s}(\R^n)$, for any $\nu\in \mathbb S^{n-1}$, we get
$$
|(-\Delta)^s w(x)|\leq |(-\Delta)^s w(R\nu)|\qquad\forall x\in \R^n\setminus\overline B_R.
$$
From \eqref{Dyda formula}, we have
$$
|(-\Delta)^s w(R\nu)|=(1+s) K_{n,s}^{-1}R^2\frac{2s}{n}.
$$
Therefore, $(-\Delta)^sw\in L^\infty(\R^n)$, and so $w$ is an admissible competitor for $\Lambda^C_{s,\infty}$. Moreover, we deduce the following estimate
\begin{align*}
\Lambda^C_{s,\infty}(B_R)&\le\frac{\left\|(-\Delta)^sw\right\|_{L^\infty(\R^n)}}{\left\|w\right\|_{L^\infty(\R^n)}}\\
&\le K^{-1}_{n,s}\max\left\{1,\frac{2s}{n}\right\}\frac{(1+s)}{R^{2s}} \\
&=(1+s)\max\left\{1,\frac{2s}{n}\right\}\Lambda^H_{s,\infty}(B_R),
\end{align*}
and this plainly proves \eqref{eq:estimatesclampedeigenvalue}.
\end{proof}

\textsc{Acknowledgements:} The authors are grateful to Serena Dipierro, Luca Gennaioli, and Enrico Valdinoci for useful discussions and suggestions. S.C. gratefully acknowledges the University of Western Australia (UWA), Serena Dipierro, and Enrico Valdinoci for their kind invitation and hospitality during the period from 24 April to 15 May 2026, when the work on this article was initiated.
S.C. acknowledges partial financial support through the EPSRC grants EP/X017206/1 and EP/X017109/1. This work was supported by the Italian Ministry of University and Research (MUR) under the FIS 2
programme, project BRIDGE: “Beyond Regularity in Differential Geometric Evolutions” (project
code FIS-2023-02224; CUP E53C25001790001). \\
The authors are members of the INdAM-GNAMPA of Italy and they have been partially supported
	by the INdAM - GNAMPA 2026 Project ``Analisi variazionale per operatori locali e nonlocali
	possibilmente singolari o degeneri''.\\

\textsc{AI statement}:  
In the paper we exploit the interactions with ChatGPT Plus to assist ourselves with stylistic improvements in the English language, LaTeX formatting, and the identification of relevant references related to the proof of the rigidity in Theorem \ref{FK hinged}. All mathematical arguments are original and are entirely the work of the authors, who are solely responsible for the mathematical content, results, and conclusions presented in this article.\\

\textsc{Data availability}: Our manuscript has no associated data.\\

\textsc{Conflict of interest}: The authors have no Conflict of interest to declare for this article.

\end{document}